\documentclass[12pt]{article}
\usepackage{amsmath,amsfonts,amssymb}
\usepackage{fullpage}

\begin{document}

\title{Siegel zeros of Eisenstein series}
\author{Joseph Hundley\footnote{the author was supported by an NSF 
VIGRE postdoctoral fellowship} 
\\ Penn State University Mathematics Department\\
University Park, State College PA, 16802, USA}
\maketitle

\newtheorem{lem}{Lemma}
\newtheorem{fact}[lem]{Fact}
\newtheorem{cor}[lem]{Corollary}
\newtheorem{thm}[lem]{ Theorem}
\newcommand{\pf}{{\bf Proof:  }}
\newcommand{\ord}{\mbox{ord}}
\newcommand{\defn}{{\bf Definition:  }}
\newtheorem{prop}[lem]{ Proposition}
\newcommand{\done}{\hfill  $\blacksquare$}
\newcommand{\donet}{ \begin{flushright}$\blacksquare$\end{flushright}}
\newcommand{\fa}{\ensuremath{\frak{a}}}\newcommand{\W}{\ensuremath{\mathcal{W}}}
\newcommand{\G}{\ensuremath{\mathbf{G}}}
\newcommand{\K}{\ensuremath{\mathbf{K}}}
\newcommand{\M}{\ensuremath{\mathbf{M}}}
\newcommand{\PP}{\ensuremath{\mathbf{P}}}
\newcommand{\C}{\ensuremath{{\mathbb C}}}
\newcommand{\A}{\ensuremath{{\mathbb A}}}
\newcommand{\uhp}{\ensuremath{{\mathbb H}}}
\newcommand{\Z}{\ensuremath{{\mathbb Z}}}
\newcommand{\Q}{\ensuremath{{\mathbb Q}}}
\newcommand{\R}{\ensuremath{{\mathbb R}}}
\newcommand{\im}{\mbox{Im}}
\newcommand{\re}{\mbox{Re}}
\newcommand{\res}{\mbox{Res}}
\newcommand{\bs}{\ensuremath{\backslash}}
\newcommand{\N}{\ensuremath{{\mathbb N}}}
\newcommand{\lamda}{\lambda}
\newcommand{\lo}{\lambda_0}
\newcommand{\lra}{\ensuremath{\longrightarrow}}
\newcommand{\wm}{\ensuremath{w_{\mbox{\tiny max}}}}
\newcommand{\Wm}{\ensuremath{W_{\mbox{\tiny max}}}}
\newcommand{\ws}{\ensuremath{W_{\mbox{\tiny sng}}}}
\newcommand{\wms}{\ensuremath{w_{\mbox{\tiny ms}}}}
\newcommand{\fm}{\ensuremath{\frak{m}}}
\newcommand{\fw}{\ensuremath{\frak{w}}}
\newcommand{\cha}{\ensuremath{\tilde{\alpha}}}
\newcommand{\I}{\ensuremath{\mathcal{I}}}

\begin{abstract}
If $E(z,s)=\sum_{(m,n)\neq(0,0)}\frac{y^s}{|mz+n|^{2s}}$ is
the nonholomorphic Eisenstein series on the upper half plane, then for
all $y$ sufficiently large, $E(z,s)$ has a "Siegel zero." That is
$E(z,\beta)=0$ for a real number $\beta$ just to the left of one.  We give
a generalization of this result to Eisenstein series formed with real
valued automorphic forms on a finite central covering of the adele points
of a connected reductive algebraic group over a global field.

Mathematics Subject Classification:    
32N05 General theory of automorphic functions of several complex variables.
\end{abstract}

\section{ Introduction}
\label{intro}
The notion of a ``Siegel zero" was first introduced in the context of Dirichlet
$L$-functions, $L(s,\chi)$, where $\chi$ is a quadratic Dirichlet character
of conductor $D.$ 
Roughly speaking, a Siegel zero of $L(s,\chi)$ is a zero
on the real axis close to 1.  More precisely, given 
$\epsilon>0$, an $\epsilon$-Siegel zero of $L(s,\chi)$ is a real number
$\beta$ in the interval 
$(1-\epsilon (\log D)^{-1},1)$  such that $L(\beta,\chi)=0$.     
The idea of a "Siegel zero" has been generalized to
$L$-functions associated with Maass forms by Hoffstein and 
Lockhart ~\cite{hofflock}.  In this context, the conductor $D$ is replaced by
$(\lambda N+1)$, where $\lambda$ and $N$ are the Laplace eigenvalue
and level of a Maass
form, respectively.

The connection with Eisenstein series comes from a 1975 paper of Goldfeld
~\cite{dg} 
in which he proves that if one assumes that the class number, $h(-d)$, of 
the imaginary quadratic field $k$ of discriminant $-d$
is small, 
then the Siegel zero is given by a certain asymptotic formula.  This 
asymptotic formula is computed by writing the Dedekind zeta function
$\zeta_k$ in 
two ways: 
\begin{equation}\label{dg}
\zeta(s)L(s,\chi_{-d})=\zeta_k(s)= \frac{1}{|\frak{o}_k^{\times}|}
\sum_C 2^s d^{-\frac{s}{2}}E(z_C,s),
\end{equation}
 where 
$E(z,s)$ is the nonholomorphic Eisenstein series on the upper half
plane, the sum
is over Heegner points $z_C$ corresponding to ideal classes $C$ of $k$,
and $\chi_{-d}$ is the unique quadratic character of conductor $d$ such that
$\chi_{-d}(-1)=-1.$
It is interesting, in this context, to 
note that for any fixed $\epsilon,$ and 
$d$ sufficiently large, a positive proportion of the functions
$E(z_C,s)$ vanish for some $s$ 
in the interval $(1-\epsilon (\log d)^{-1},1)$.  
To prove this one uses the equidistribution of Heegner points 
and the following result:

\begin{thm} ($GL(2)$ Case, due to Bateman and Grosswald \cite{bg}):  
Given $\epsilon >0,$ there exists $Y >0$ such that
if $y>Y$ then for each $x,$ $E(x+iy, \beta(x,y)) =0$ for some 
$\beta(x,y) \in (1-\epsilon(\log y)^{-1},1),$ and for all $x,\epsilon$,
we have
$$1-\beta(x,y)\sim \frac{3}{\pi y}, \mbox{ as } y \rightarrow \infty.$$ 
\end{thm}

Hoffstein \cite{jhoff} has proved an analogous result for Hilbert modular
Eisenstein series.  We extend it to the generality of
Moeglin and Waldspurger \cite{mw}.  However, our 
result will be weaker than Hoffstein's in that we will not obtain a 
precise error term, and in that we will not obtain any estimates for the 
optimal value of $Y$.  

The author wishes to thank Carlos Moreno, who suggested that the result, 
originally proved only for $GL_n$ over $\Q$ would go through in this 
generality, as well as Dorian Goldfeld, Herv\'e Jacquet, and David 
Ginzburg, for help and advice along the way.

\subsection{Notation}
\label{notation}

For the most part, we follow the notation of Moeglin and
Walspurger \cite{mw}.  This entails 
certain redundancies:  for example the symbol $k$ is used both for
the global field and an element of the maximal compact, while $M$ is 
both a Levi subgroup and an intertwining operator.  All references
are to \cite{mw} until section 4

Thus, let $k$ be a global field, $\A$ the adeles of $k$,
$G$ a
connected reductive algebraic group defined over $k$, and $\mathbf{G}$
a finite central covering of $G(\A)$ such that $G(k)$ lifts to a 
subgroup of $\G$.

Fix once and for all a choice of minimal parabolic $P_0$ of $G$ 
and a Levi subgroup
$M_0$ of $P_0,$ both defined over $k$.    
This also fixes a definition of ``standard'' for parabolics and Levis (p.4).
Let $T_0$ be the maximal split torus in the center of $M_0$ and
$\Delta_0$ the set of simple positive roots for $T_0$ determined by
$P_0.$  
Fix a maximal compact subgroup $\mathbf{K}$ of $\G$,
as in I.1.4.

Let $P$ be a standard parabolic subgroup of $G$, and let $U$ be its 
unipotent radical, and $M$ its standard Levi. 
Then $U(\A)$ lifts canonically into $\G$.  
Let $\M$ denote the preimage of $M(\A)$ in $\G$.  
If $\chi$ is a rational character of $M$, we obtain a map 
$\M\mapsto \R^{\times}_+$ by projecting to $M(\A)$, using $\chi$ to 
get to $\A^{\times}$, and then taking the absolute value.  
Let $\M^1$ denote the intersection of the kernels of all such maps.
Let $X_M^{\G}$ and $\re X_M^{\G}$ 
denote the groups of continuous homomorphisms of $\M$ to
$\C^{\times}$ and $\R^{\times}_+$ respectively, which
are trivial on $\M^1$ and on the center of $\G$.  When $k$ is a 
number field, $X_M^{\G}$ may be identified with a complex vector
space $(\frak{a}_M^{\G})^*$.  When $k$ is a function field, it has 
a finite number of connected components, and there is a natural 
projection from $(\frak{a}_M^{\G})^*$ to the identity
component, which still identifies $\re X_M^\G$ with a real vector 
space.  See I.1.4 and I.1.6.

There is unique map $m_P:\G\longrightarrow \M^1\bs \M$ defined by requiring that if $g=umk$ for
some $u\in U(\A), m\in \M$ and $k \in \K$, then $m_P(g) =\M^1m.$  
If $\varphi$ is a function $U(\A)M(k)\bs \G \mapsto \C$ and $k\in \K$, 
let us define a function $\varphi_{k}$ on $M(k)\bs \M$ by 
$\varphi_{k}(m)=m^{-\rho_P}\varphi(mk),$ where $\rho_P$ is half the sum of the
roots of $M$ in Lie $U$.  
Let $\phi_{\pi}$ be an automorphic form on
$U(\A)M(k)\bs \G$ such that
for each $k$, the function $\phi_{\pi,k}$ is a cusp form on $\M$ which 
generates a
semisimple isotypic submodule of type
 $\pi,$  where $\pi$ is an automorphic subrepresentation of 
$\M$ in the sense of \cite{mw}, page 78.

For each 
$\lambda \in X_M^{\G}$, let 
$$\lambda\phi_{\pi}(g)=m_{P}(g)^{\lambda}\phi_{\pi}(g).$$  Then for each $k$, the function $\lambda \phi_{\pi,k}$ is a cusp form on $\M$ which generates a semisimple isotypic submodule of type $\pi\otimes \lambda.$

For $\lambda$ in a suitable cone in $X_M^{\G}$, the Eisenstein series is defined by the following convergent sum:
$$E(\lambda \phi_{\pi},\pi\otimes \lambda)(g)=\sum_{\gamma \in P(k)\bs G(k)}\lambda \phi_{\pi}(\gamma g).$$  It is holomorphic for $\lambda$ in the domain of convergence, and extends to a meromorphic 
function on all of $X_M^{\G}$.  We may assume without loss of generality
that $\pi$ is unitary.  This amounts to,
at most, altering our choice of ``base point.''  (See  I.3.3)
(We deviate from \cite{mw} in viewing the 
Eisenstein series as a function on $X_M^{\G}$, rather than their $\frak{P}$,
which is a principal homogenous space for the quotient of $X_M^{\G}$ by 
a certain finite group.) 

Suppose that $\phi_{\pi}$ is real valued.  Then 
$E(\lambda \phi_{\pi},\pi\otimes \lambda)$ is real valued for 
$\lambda \in \re X_M^{\G}.$
Let us fix $P$ and $\phi_{\pi}$ once and for all.  
For each $\alpha \in \Delta_0,$ let 
$P_{\alpha}=M_{\alpha}U_{\alpha}$
denote the standard maximal parabolic such that every element of 
$\Delta_0$ is a root of $M_{\alpha}$ except $\alpha.$ 
Our choice of $\Delta_0$ determines a set of positive roots for 
the action of $T_0$ on any Levi $M'$ which we denote by $R^+(T_0,M')$.

Let $W$ denote the Weyl group of $G$, and $W_M$ that of $M$.  
Let $W(M)$ be the set of $w\in W$ of
minimal length in their class $wW_M$, such that $wMw^{-1}$ is a 
standard Levi of $G$.
Let $W(M,M_{\alpha})$ denote the set of $w \in W$ such that $w^{-1}\theta >0$ for every $\theta \in R^+(T_0,M_{\alpha})$ and $wMw^{-1}$ is a standard 
Levi subgroup of $M_{\alpha}.$ 

If $k$ is a function field we fix once and for all a place $v_0$, 
and a uniformizing parameter $\frak{w}$, and let $q=|\frak{w}|^{-1}_{v_0}.$  
Let
$\frak{m}=\R_+^{\times}$ in the number
field case or $\frak{w}^{\Z}$, in the function field case.  Thus $\fm$ may be 
embedded in $\A^{\times}$ either at $k_{v_0}$ or diagonally at the infinite 
places, as a subgroup on which  
 the absolute value is injective.  One then has a 
subgroup of $T_0(\A)$ isomorphic to $\fm^R$, where $R$ is the rank of $T_0$, 
and
in  I.2.1 of \cite{mw} this is extended to a subgroup $A_{\M_0}$ of $\G$, 
still isomorphic to $\fm^R.$  We then define $A_{\M_{\alpha}}=A_{\M_0}\cap 
Z_{\M_{\alpha}}.$  

If $\lambda \in \re X_M^{\G}$, then $\left. \lambda \right|_{A_{\M_{\alpha}}}$
factors through the quotient $A_{\M_{\alpha}}/A_\G$, and is trivial on any 
torsion in this quotient.  It follows that 
we may fix a map
$\tilde{\alpha}:\fm \rightarrow A_{\M_{\alpha}}$ such that
$\left. \lambda \right|_{A_{\M_{\alpha}}}$ is determined by $\lambda \circ 
\tilde{\alpha}.$  
In the function field case, we also denote by $\tilde{\alpha}$ the 
composition of $\tilde{\alpha}$ with  $|\cdot |_{v_0}^{-1}:q^{\Z}
\rightarrow \frak{w}^\Z.$
Then $\lambda \circ \cha$ is a continuous homomorphism of a subgroup of 
$\R_+^{\times}$ to $\R_+^{\times}$, so 
we may define $\langle\tilde{\alpha},\lambda\rangle$ to be the 
unique real number such that
$$\lambda\circ\tilde{\alpha}(y)=y^{\langle\tilde{\alpha},\lambda\rangle}
.$$
Replacing $\tilde{\alpha}$ by the map $y\mapsto \tilde{\alpha}(y^{-1})$ 
if necessary, we may add the stipulation that 
$\langle \tilde{\alpha},\alpha\rangle >0.$
We will also need to refer to the Siegel set $S$, defined in section I.2.1 of 
\cite{mw}, using a compact subset $\omega$ of $\PP_0.$

\defn Let us say that a map $\Lambda : \C\rightarrow X_M^{\G}$ is 
\emph{elementary} if it is holomorphic and the restriction to $\R$ is an 
affine map into the real vector space $\re X_M^\G.$

The constant term of $E(\lambda\phi_\pi,\lambda\otimes\pi)$ along $P_\alpha$
(see I.2.6)
is given in terms of Eisenstein series $E^{M_{\alpha}},$ 
defined analogously, with $M_{\alpha}$ replacing $G$, and
intertwining operators $M(w,\pi)$ defined in  II.1.6:
\begin{equation}\label{ctfrommw}
E_{P_{\alpha}}(\lambda\phi_{\pi},\lambda\otimes\pi)=
\sum_{w\in W(M,M_{\alpha})}E^{M_{\alpha}}
(M(w,\lambda\otimes \pi)\lambda\phi_{\pi},w(\lambda \otimes \pi)).
\end{equation}
(See II.1.7).

Suppose that $E$ has a singularity along a root 
hyperplane $H$, associated to a root $\theta$ as in IV.1.6.  
Our theorem applies only to the case when $H_\R:=H\cap \re X_M^\G$ 
is non-empty.
In this case, $H_\R$ will be of the form
$$\{\lambda \in \re X_M^{\G}: \langle \lambda, \check{\theta} \rangle=c_H\},$$ 
for some $c_H\in \R$.  See
IV.1.6 and I.1.11.  
Here $\check{\theta}$ may be taken to be the coroot associated to a positive
root $\theta,$ or its projection to the dual of $\re X_M^\G$:  the set $H$ is
the same either way.
We say that $\lambda\in H$ is \emph{generic for $\alpha$} 
if $\lambda$ does not lie in 
any other hyperplane along which $E$ has a singularity, and
$\langle\tilde{\alpha},w_1\lambda\rangle = 
\langle\tilde{\alpha},w_2\lambda\rangle$ for  $w_1,w_2\in W(M,M_{\alpha})$
iff $w_1=w_2$.
It follows that for $W(M,M_\alpha)$ nonempty, and $\lambda$ generic 
for $\alpha,$
there is a unique 
$\wm(\lambda,\alpha)$ 
such that $\langle\tilde{\alpha},w\lambda\rangle$ is maximal. 

\subsection{Statement of Main Result}
\label{statement}

\begin{thm}\label{mainthm}
Fix a root hyperplane $H$ along which the Eisenstein series is singular, 
a root $\alpha$,
an elementary map $\Lambda$ such that $\Lambda(\C)\cap H=\Lambda(0)$ is
generic for $\alpha$,  and elements $p\in \omega, k\in \K.$ 
Let $$E(y,s)=E(\Lambda(s)\phi_{\pi}, \Lambda(s)\otimes \pi)
(p\tilde{\alpha}(y)k).$$
and for each $w\in W(M,M_\alpha),$ let 
$$E_w(y,s)=E^{M_\alpha}(
M(w,\Lambda(s)\otimes \pi)\Lambda(s)\phi_{\pi},w( \Lambda(s)\otimes \pi))
(p\tilde{\alpha}(y)k).$$
We say that $\beta$ is an $\epsilon$-Siegel zero of $E(y,\sigma)$ if
$\beta \in (-\epsilon (\log y)^{-1},\epsilon (\log y)^{-1})$ and
$E(y,\beta)=0.$
Let $\ws(H,\alpha,\Lambda,pk)=\{w\in W(M,M\alpha): E_w(1,s) \mbox{ has a 
pole at }s=0\}.$
If $H,\alpha,p,k,$ and $\Lambda$ satisfy
\begin{enumerate}
\item{The set $W(M,M_\alpha)$ is nonempty,
and $\wm(\Lambda(0),\alpha)\notin \ws(H,\alpha,\Lambda,pk),$}
\item{$E_{\wm(\Lambda(0),\alpha)}(1,0)\neq 0,$}
\item{$E(y,s)$ has a simple pole at $s=0.$}
\end{enumerate}
then we have the following
\begin{description}
\item[A)]{ 
If $\ws(H,\alpha,\Lambda,pk)\neq \emptyset$, 
let
$\wms$  be the unique element of $\ws$ such that 
$\langle \tilde{\alpha},\wms\Lambda(0) \rangle$ is maximal 
(among elements of \ws), and denote $\wm(\Lambda(0),\alpha)$
more briefly by $\wm.$  Then
\begin{description}
\item[i)]{There exists $Y>0$ 
(dependent on $\epsilon,H,\alpha,p,k,\Lambda,$ and $\phi_\pi$), 
such that for all $y>Y$, $E(y,\sigma)$
has an $\epsilon$-Siegel zero.
}
\item[ii)]{Let $\beta:(Y,\infty)\rightarrow \R$ or $q^\Z\cap
(Y,\infty)\rightarrow \R$ in the function field case,
be any function such that for
each $y$, $\beta(y)$ is an $\epsilon$-Siegel zero of $E(y,\sigma).$  
Then
$\beta(y) \sim e y^{-\langle \cha, \wm \Lambda(0)-\wms \Lambda(0) \rangle},$
where 
$$e=
 \frac{\chi_\pi(\wms^{-1}\cha(y)\wms)
\left({{\res}\atop{s =0}}E^{M_\alpha}(M(\wms,\Lambda(s)\otimes \pi)\Lambda(s)
\phi_\pi,\wms(\Lambda(s)\otimes \pi))(pk)\right)
}{\chi_\pi(\wm^{-1}\cha(y)\wm)
E^{M_\alpha}(M(\wm,\Lambda(0)\otimes \pi)\Lambda(0)
\phi_\pi,\wm(\Lambda(0)\otimes \pi))(pk)
}.
$$
}
\end{description}}
\item[B)]{
If $\ws(H,\alpha,\Lambda,pk)=\emptyset$,
then the conclusion is 
weaker:
if
$y>Y$ then either $E(y,\sigma)$ has an $\epsilon$-Siegel zero, or
$E(y,\sigma)$ is holomorphic at zero.
If $k$ is a function field, then it is always the latter.}
\end{description}
\end{thm}
We will also prove the following lemma, relevant to the problem 
of choosing $\alpha,\Lambda,p$ and $k$ so that conditions 1-3 above 
are satisfied.  
\begin{lem}\label{condlem}
For $\alpha \in \Delta_0$ such that $\alpha \notin R^+(T_0,M)$ and 
$\theta\notin R^+(T_0,M_\alpha)$, we have
\begin{enumerate}
\item{the identity element, $1$, is in $W(M,M_\alpha),$}
\item{for all $p,k$ and $\Lambda$, $1 \notin \ws(H,\alpha,\Lambda,pk),$}
\item{if $c_H>0$, then 
$\{\lambda \in H_\R,\mbox{ generic for }\alpha: \wm(\lambda,\alpha)=1\}$ is a 
nonempty open subset of $H_\R$.}
\end{enumerate} 
\end{lem}
It is also known (IV.1.11(c)) 
that when $c_H>0,$ the singularity along $H$ is without
multiplicity.  For each $p,k$ the restriction of $E_1$ and the 
continuation of $h_H(\lambda)E$ are meromorphic functions of $\lambda
\in H$.  
As we see in an example below, either or both may be trivial 
(i.e., zero for all $\lambda$) for 
certain values of $p,k$.  But, for any $p,k$ such that both are 
nontrivial, $\Lambda$ may
then be chosen so that all three conditions are satisfied.  

\section{ Basic Analytic Result}
\label{bar}
Our approach to proving the existence of Siegel zeros of Eisenstein series is
as follows.  First, we show that any function that can be written in a 
certain form will have a zero on the real axis very close to its pole. Then, 
we show that our Eisenstein series can always be put into that form.  
The lengthy definition of the function $F(y,s)$ in the following lemma should,
therefore, be taken as the ``mold'' into which we will fit the Eisenstein 
series.  As we see here, any function which fits into this mold will have
Siegel zeros.\label{secwitlem}

\begin{lem}\label{mainlemma}
Fix real numbers $a, b, c,$ and $d,$ such that 
$b > d$.  
Let $\sigma$ be a real variable, 
which we may think of as restricted to a small neighborhood of $0$ 
and
let $y$ be another real variable, which we think of as positive and large,
keeping in mind that when we apply this Lemma to the function field case, 
$y$ will only range over $q^{\Z}.$  Let 
$A(\sigma)$ and $C(\sigma)$ be two real valued 
functions that are both continuous and nonvanishing
for $\sigma$ in a neighborhood of $0,$ 
and 
let  $B(y,\sigma)$ and $D(y,\sigma)$ be two more real valued functions, 
which are 
both continuous in $\sigma$, and such that $B(y,\sigma)y^{-(a\sigma+b)}$, and 
$D(y,\sigma)y^{-(c\sigma+d)}$ tend to zero as $y$ tends to $\infty$ for all values 
of $\sigma$ in some 
neighborhood of $0$, and that convergence is uniform as $\sigma$ ranges
over this neighborhood.  

Define:
$$F(y,\sigma)=(A(\sigma)y^{a\sigma+b}+B(y,\sigma))+\frac{1}{\sigma}
(y^{c\sigma+d}C(\sigma)+D(y,\sigma)).$$

Then we have
\begin{description}\item[i)]{
 For every $\epsilon>0,$ there exists $Y(\epsilon)>0,$ 
such that if $y>Y(\epsilon),$ 
then $F(y,\sigma)$ has a zero 
in the interval $(-\epsilon (\log y)^{-1},\epsilon (\log y)^{-1})$.}
\item[ii)]{
Now fix an $\epsilon >0$, and 
take $\beta: (Y(\epsilon),\infty)\rightarrow\R$ such that for each $y$ 
we have 
 $\beta(y)\in (-\epsilon (\log y)^{-1},\epsilon (\log y)^{-1}),$ and
$F(\beta(y),y)=0$.  For any such $\beta$, 
$$ 
-\beta(y) \sim \frac{C(0)}{y^{(b-d)}A(0)} \mbox { as }
y \rightarrow \infty.$$} \end{description} 
\end{lem}

\pf  
By replacing 
$F(y,\sigma)$ with 
$y^{-(c\sigma+d)}F(y,\sigma)$, which has the same zeros, we may assume $c=d=0.$
By considering the four functions $\pm F(y,\pm \sigma)$ we may assume that
$A(0)$ and $C(0)$ are positive.  

Next, choose $\delta, Y_1, m, M$ such that 
\begin{equation} \label{yano} 
y>Y_1, |\sigma|<\delta \Rightarrow \left\{\begin{array}{l}
0<m<A(\sigma)<M,\\0<m<C(\sigma)<M, \\ |D(y,\sigma)|<\frac{m}{2}, \\ 
|B(y,\sigma)|<\frac{m}{2} y^{a\sigma+b}.
\end{array}\right.
\end{equation}
Then, for each $y>Y_1,$ the function
$A(\sigma)y^{a\sigma+b}+B(y,\sigma)$ is bounded on 
$|\sigma|<\delta,$
while $C(\sigma)+D(y,\sigma)$ is bounded away from zero.
Hence, for
every such $y$, there is a neighborhood of the form $(-l,0)$ on which 
$F(y,\sigma)<0$.
Now suppose that $y>Y_1$, and 
$\epsilon (\log y)^{-1} <\delta.$  Then the bounds (\ref{yano}) are valid at
$\sigma=-\epsilon (\log y)^{-1},$ yielding
$$F(y,-\epsilon (\log y)^{-1})>\frac{m}{2}e^{-a\epsilon} y^{b}
-\epsilon^{-1} (M+\frac{m}{2}) \log y.$$ 
Clearly, if we fix an $\epsilon >0$, we may first choose $Y_2\geq Y_1$ such
that  
the above is valid for all $y>Y_2$, and then
 choose $Y_3$ such that if $y>Y_3,$ the right side is
positive.  We then let $Y=\max(Y_2,Y_3)$, 
and the first assertion is proved. 
Moreover, if we fix and $\epsilon'<\epsilon$, we can choose $Y_2'$ 
such that for every $y>Y_2', \epsilon'<\varepsilon<\epsilon,$ $\varepsilon 
(\log y)^{-1} <\delta,$ and then choose $Y_3'$ such that for
$y>Y_3', \epsilon'<\varepsilon<\epsilon,$ 
$$\frac{m}{2}e^{-a\varepsilon} y^{b}
-\varepsilon^{-1} (M+\frac{m}{2}) \log y >0.$$  
It follows that any $\beta(y),$ defined as above \emph{with respect to}
$\epsilon$ satisfies
$$\beta(y) \in (-\epsilon'(\log y)^{-1},0) \mbox { for }y>Y'.$$
Since this works for any $\epsilon',$ we have
$$\lim_{y \rightarrow \infty} y^{-\beta(y)}=1$$
independently of the choice of $\epsilon,$ and independently of any possible 
choice of $\beta(y).$
To prove \emph{ii)}, we note that 
$F(y,\beta(y))=0$ iff  $$A(\beta(y))y^{a \beta(y) +b} +B(y,\beta(y))=
\frac{-1}{\beta(y)}(C(\beta(y))+D(y,\beta(y))).$$  We have seen that 
for $\beta(y)$ near $0$ and $y$ large the left side is nonzero, so we may 
put this into the form
$$\beta(y)=-\frac{C(\beta(y))+D(y,\beta(y))}{A(\beta(y))y^{a \beta(y) +b} +B(y,\beta(y))}.$$
So
\begin{eqnarray*}&&
\hskip -.2in \lim_{y\rightarrow \infty}
\frac{C(0)/(A(0)y^{b})}{-\beta(y)}
\\&&=\lim_{y\rightarrow \infty} \frac{C(0)}{C(\beta(y))+D(y,\beta(y))}
\frac{\left(A(\beta(y))y^{a\beta(y)+b}+B(y,\beta(y))\right)y^{-b}}{A(0)}.
\end{eqnarray*}
The second limit is evidently 1, which proves the asymptotic formula.  \donet

{\bf Remarks 1.} If, on the other hand, $F$ is in the same form, with 
$b < d,$ but all
the other assumptions are the same, then similar arguments show that 
$F(y,\sigma)$ will be nonvanishing 
for $\sigma$ in a neighborhood of $0$ and $y$ sufficiently large.

{\bf 2.} If $F$ is as above, but $C(0)=0,$ then the asymptotic
formula is no longer correct, but the Siegel zero still exists
for all values of $y>Y$ such that $D(y,\sigma)\neq 0.$ 

\subsection{Necessary Fact}
\label{necessaryfacts}
We will need one more well known fact from the theory of Eisenstein series: 
the Eisenstein series is well approximated by its constant term.  
\begin{prop}
\label{approx}
If $k$ is a function field, then there is a constant $c$ such that
$E(\lambda\phi_{\pi},\pi\otimes\lambda)(g)-
E_{P_{\alpha}}(\lambda\phi_{\pi},\pi\otimes\lambda)(g)=0$,
whenever $g\in S$ and $m_{P_0}(g)^{\alpha}>c$, where $S$ is a Siegel domain
as in I.2.1.  
In particular, 
for $\Lambda,p,k,\cha$ as in the main theorem,
$$(E-E_{P_{\alpha}})(\Lambda(\sigma)\phi_{\pi},\pi\otimes\Lambda(\sigma))
(p\tilde{\alpha}(y)k)=0,$$
for $p\in \omega, k\in K$ and $y$ sufficiently large.
If $k$ is a number field, then for $\Lambda,p,k,\cha$ as in the main theorem,
$$
\sigma\left(
(E-E_{P_{\alpha}})(\Lambda(\sigma)\phi_{\pi},\Lambda(\sigma)
\otimes\pi)(p\tilde{\alpha}(y)k)
\right)
$$ 
is rapidly decreasing as a function of $y$ for $p\in \G^1\cap \omega,k\in K$, 
uniformly for $\sigma$ in a neighborhood of zero.
\end{prop}
\pf
The function field case is immediate from Lemma I.2.7.  
For the number field case we use Lemma I.2.10.  
Lemma I.4.4, in conjunction with I.2.5 yields the bounds required by
the hypotheses of Lemma I.2.10.
\done
\section{Proofs}
\subsubsection{Proof of Theorem \ref{mainthm}}
As $H,\alpha,\Lambda,p$ and $k$ are fixed, we suppress them from the
notation, denoting $\ws(H,\alpha,\Lambda,pk)$ by $\ws$, etc.
The idea is to fit $E(y,\sigma)$ into the mold described in section 3.
In the new notation we have introduced, equation \eqref{ctfrommw} reads
\begin{equation}\label{decomp}
E_{P_{\alpha}}(y,\sigma)=\sum_{w\in W(M,M_{\alpha})}E_w(y,\sigma).
\end{equation}
We observe that
\begin{equation}\label{athing}E_w(y,\sigma)=
\chi_\pi(w^{-1}\tilde{\alpha}(y)w)
y^{\langle\tilde{\alpha}, \rho_{P_\alpha}+w \Lambda(\sigma)\rangle}
E_w(1,\sigma).
\end{equation}
Now, $\chi_\pi$ is both real-valued and unitary.  When $k$ is a number
field, 
$\chi_\pi(w^{-1}\tilde{\alpha}(y)w)$
 is trivial, but when $k$ is a function field, 
$\chi_\pi(w^{-1}\tilde{\alpha}(q^n)w)$ may equal $-1$ for some $w$ when $n$ is
odd.  In this case, one may consider restrictions to odd and even $n$ 
seperately, and combine the results.  We omit the details, and assume that
$\chi_\pi(w^{-1}\tilde{\alpha}(y)w)$ is identically $1$, when $w$ is either
$\wm$ or $\wms.$

For each $w$, the map 
$\sigma \mapsto \langle \tilde{\alpha},w \Lambda(\sigma) \rangle$ 
is an affine map $\R \rightarrow \R$, and so we may define 
$a,b,c,$ and $d$ by the conditions 
$\langle \tilde{\alpha}, \rho_{P_\alpha}+\wm\Lambda(\sigma)\rangle
=a\sigma+b,$ and,
if $\ws$ is nonempty,
$\langle \tilde{\alpha}, \rho_{P_\alpha}+\wms\Lambda(\sigma)\rangle
=c\sigma+d.$
If $\ws$ is empty, we take $c=0$, and $d$ any real number less than 
$b$.
We take 
\begin{eqnarray*}
A(\sigma)&=&E_{\wm}(1,\sigma),\\
C(\sigma)&=&\sigma E_{w_{\mbox{\tiny ms}}}(1,\sigma), \mbox{ or }
0 \mbox{ if }\ws=\emptyset,\\
D(y,\sigma)&=&\sigma \left(
\sum_{w \in \ws-\left\{w_{\mbox{\tiny ms}}\right\}}E_w(y,\sigma)
+(E-E_{P_\alpha})(y,\sigma)\right)
,\\
B(y,\sigma)&=&\sum_{w \in W(M,M_{\alpha})-\ws-\left\{\wm\right\}}E_w(y,\sigma),
\end{eqnarray*}
Where we extend $C$ and $D$ to $\sigma=0$ by continuity.  

We check that all the hypotheses of the main lemma are satisfied.  
As, $\Lambda(0)$ is generic, 
we may choose an open ball containing it 
which intersects no other hyperplane along which $E$ is singular.  
The continuity of each function on the set of $\sigma$ mapping into 
this ball is clear.  The fact that all functions are real valued follows
from the fact that $\phi_\pi$ is.  Condition 2 is precisely that 
$A(0)\neq 0,$ and in the case when $\ws \neq \emptyset,$ the value 
of $C$ at $0$ is the residue of $E_{w_{\mbox{\tiny ms}}}(1,\sigma)$
at zero.  The bounds on $B$ and $D$ as $y \rightarrow \infty$ follow
easily from \eqref{athing} and Proposition \ref{approx}.\done

\subsubsection{Proof of Lemma \ref{condlem}}
The first assertion is trivial.
The second and third follow from more general assertions, which we now
prove.
For $w \in W(M,M_\alpha)$, let $W_{M_\alpha}(wMw^{-1})$ be defined in the same
way as $W(M),$ with $M_\alpha$ replacing $G$ and $wMw^{-1}$ replacing $M$.

\begin{lem} \label{hollemma}
 Let $H$ be a root hyperplane associated to a root $\theta$,
and fix $\alpha.$
If there exist $p,k,\Lambda$ such that 
$\Lambda(0)\in H$ is generic for $\alpha$, and
$w \in \ws(H,\alpha,p,k,\Lambda)$, 
then there exists $w'\in W_{M_\alpha}(wMw^{-1})$ such that 
$w'w\theta<0.$
\end{lem}

\pf  For some
$w'\in W_{M_\alpha}(wMw^{-1}),$ write $w'w=s_{\gamma_\ell} \dots
s_{\gamma_1}$, where for each $i$, $s_{\gamma_i}$ is an ``elementary
symmetry'' as in section I.1.8.  See also section IV.4.1.  Let
$w(j)=s_{\gamma_j}\dots s_{\gamma_1}.$  Then
$$M(w'w,\pi \otimes \lambda)=M(s_{\gamma_\ell},w(\ell-1)(\pi \otimes \lambda))
\circ \dots \circ M(s_{\gamma_1},\pi\otimes\lambda).$$
The singularities of $M(s_\gamma,\tau \otimes \mu)$ are carried by a 
locally finite set of hyperplanes associated to $\gamma.$  Hence the 
singularities of $M(s_{\gamma_i},w(i-1)(\pi\otimes\lamda))$ are carried by
a locally finite set of root hyperplanes associated to the root 
$w(i-1)^{-1}\gamma_i,$ and, in general, the singularities of 
$M(w'w,\pi\otimes\lambda)$ are carried by a locally finite set of root 
hyperplanes each of which is associated to a root $\gamma$ such that 
$w'w\gamma <0.$  

Suppose that $w'w\theta >0$ for all $w'\in W_{M_\alpha}(wMw^{-1}).$  
Then there is a locally finite set of root hyperplanes, not containing
$H$, that carries the singularities of $M(ww',\pi \otimes \lambda)$ 
for every $w'.$  Hence it carries the singularities of all the 
cuspidal components of all the constant terms of 
$E^{M_\alpha}(M(w,\pi\otimes \lambda)\lambda \phi_\pi,w(\pi\otimes \lambda))$, 
(equation \eqref{ctfrommw} above, and  IV.1.9 (b)). By I.4.10 it
carries the singularities of
$E^{M_\alpha}(M(w,\pi\otimes \lambda)\lambda \phi_\pi,w(\pi\otimes \lambda))$
as well.  The result follows.
\done

Under the hypotheses of Lemma \ref{condlem}, 
$w'\theta >0$ for all $w'\in W_{M_\alpha}(M),$ and
2 follows.  As for 
3, it's clear that the set in question is always open, and that it's 
empty iff $\{\lambda \in H_\R:\wm(\lambda,\alpha)=1\}$ is.
Thus, we only need to prove that this latter set is nonempty.

\begin{lem}
If $c_H>0$, and $1\in W(M,M_\alpha)$, then there exists $\lambda \in H_\R$
such that $\wm(\lambda,\alpha)=1.$
\end{lem}
\pf The linear functional on
$\Delta_0$ given by pairing with $\tilde{\alpha}$ corresponds to a 
positive multiple of the one given by pairing with the
fundamental
coweight $\hat{\varpi}_{\alpha}$, in the basis dual to $\Delta_0.$
Given an element $\lambda$ of $\re X_{M_0}^\G,$ which is 
irrational, i.e., has trivial kernel in the coweight lattice,
we can define
notions of $\lambda$-positivity for coroots and $\lambda$-dominance for
coweights.  If the definitions of positivity for coroots that come from
$\lambda$ and $\Delta_0$ coincide, then the definitions of dominance 
for coweights will as well.  Hence, whenever
$$\langle \lambda ,\check{\gamma} \rangle >0 \hskip .5in 
 \forall \gamma \in \Delta_0,$$
and $\lambda$ is irrational, 
we have
$$\langle \lambda, \tilde{\alpha} \rangle \geq 
\langle w\lambda ,\tilde{\alpha} \rangle 
\hskip .5in \forall  w \in W,$$
with equality only if $w^{-1}\cha=\cha.$
Now project $\lambda$ to $\bar{\lambda}\in \re X_M^\G.$  This projection 
corresponds to restriction from $Z_{\M_0}$ to $Z_{\M}$ (see pp.7,11).  
Since the image of $w^{-1}\tilde{\alpha}$ is contained in $Z_\M$ whenever 
$w\in W(M,M_\alpha)$, we find that 
$$\langle \bar{\lambda},\tilde{\alpha}\rangle \geq 
\langle w\bar{\lambda},\tilde{\alpha}\rangle
\hskip .5in \forall w\in W(M,M_\alpha),$$ 
with equality only if $w=1.$  
Thus, for $\bar{\lambda}$ the projection of any $\lambda$ as above, we have
$\wm(\bar{\lambda},\alpha)=1.$  Now, let $\gamma_\theta$ be the root 
associated to $\theta$ as in section I.1.11 (so that $\check{\theta}$ is 
defined as the projection of $\check{\gamma_\theta}$ to $\re \frak{a}_M$).
It should be emphasized that $\langle \lambda,\check{\gamma_\theta} \rangle$
and $\langle \bar{\lambda},\check{\theta} \rangle$ need not be equal.  
However, $\langle \lambda, \check{\theta} \rangle$ and 
$\langle \bar{\lambda}, \check{\theta} \rangle$ \emph{are} always equal.  
So, it's enough to verify that for some irrational $\lambda$, 
with $\langle \lambda, \check{\gamma} \rangle >0 \forall \gamma \in \Delta_0,$
the quantity
$\langle \lambda, \check{\theta} \rangle$ is also positive.  This is
straightforward, when $\check{\theta}$ is written in terms of the basis
$\Delta_0^\vee$ of coroots $\check{\gamma}$,
and $\lambda$ in terms of the dual basis.
\done

\section{Examples}
\label{examples}

Let us consider the example when $G=GL(n)$, the representation $\pi$ is unramified, and we choose $\phi_{\pi}$ to be the spherical vector.  In this case, 
the intertwining operators may be given explicitly in terms of automorphic $L$-functions, as in \cite{ep}.  It follows that a pole of the Eisenstein series comes from either a zero or a pole of one of these $L$-functions.  
The zeroes are quite mysterious, but the poles of 
Rankin-Selberg $L$-functions on $GL(n)$ have been completely determined by 
Jacquet and Shalika \cite{jsi},\cite{jsii}, and referring back to the 
$GL(2)$ case, we see that it is these poles which provide the correct 
analog of the pole of the $GL(2)$ Eisenstein series at 1. 

As noted, the case $G=GL(2)$ is the simplest example of our theorem.  Let us consider two more which exhibit some of the features of the general case, but are not so complicated as to become unwieldy.  
\subsection{Example One:  $GL(4)$ over $\Q$}

First let us consider the case $k=\Q,$ and $\G=GL(4,\A)$.  We take the 
usual Borel, consisting of all upper triangular elements, and let $P$ 
be the standard maximal
parabolic such that $M\cong GL(2)\times GL(2).$

In this case, $\re X_M^{G}$ is one dimensional, consisting of all real powers
of the modulus
$$\left(\begin{smallmatrix} h_1 &*\\ & h_2\end{smallmatrix}\right)
\mapsto \left|\frac{\det h_1}{\det h_2}\right|^2,$$
so we identify it with $\R$.  In order to have good knowledge of the poles of the intertwining operators, we will want to make a suitable choice of $\pi$ and $\phi_{\pi}$, which can also be described quite explicitly.

Let $\varphi$ be a real-valued Maass Hecke eigenform of level 1,
right invariant by the center and the maximal compact, and define a 
function on $GL(4,\A)$ by 
$$I\left(\left(\begin{smallmatrix}h_1 &*\\&h_2\end{smallmatrix}\right) k
,s,\varphi\right)=
\left|\frac{\det h_1}{\det h_2}\right|^{2s+1} \varphi(h_1)\varphi(h_2),$$
for $k$ now in the maximal compact of $GL(4,\A).$
Our Eisenstein series is 
$$E(g,s,\varphi)= \sum_{\gamma \in P(\Q)\bs
GL(4,\Q)}I(\gamma g, s,\varphi).$$
In this case $W(M,M')=\emptyset$, unless $M'=M,$ and
$$E_P(g,s,\varphi)=I(g,s,\varphi)+
\frac{L(2s,\varphi\times\varphi)}{L(2s+1,\varphi\times\varphi)}
I(g,-s,\varphi).$$
Here $L$ denotes the completed $L$ function-- including gamma factors.
(This is essentially a special case of the computation in \cite{ep}.)
Since $\re X_M^G$ is one dimensional, a codimension 1 subspace is a point.  
Thus our ``root hyperplane'' is just $s=\frac{1}{2},$ or some zero of 
$L(2s+1,\varphi\times \varphi).$  We consider the one at $\frac{1}{2},$  
and let
$\Lambda(\sigma)=\sigma+\frac{1}{2},$ and $\cha(y)=diag(y,y,1,1),$ embedded
at the infinite place.  Note that if $p=\left(\begin{smallmatrix} h_1&*\\
&h_2\end{smallmatrix}\right)$ with $\varphi(h_1)$ or $\varphi(h_2)=0,$
then the constant term is zero for all $s$.  
For any other $p$,
$\ws=
\left\{\left(\begin{smallmatrix}&I\\I&\end{smallmatrix}\right)\right\}.$
For $\sigma$ positive, $\wm=I,$  and the theorem applies.

Let $A(\sigma)=I(pk,\sigma+\frac{1}{2},\varphi),
C(\sigma)=\frac{L(2\sigma+1,\varphi\times\varphi)}
{L(2\sigma+2,\varphi\times\varphi)}I(pk,-\sigma-\frac{1}{2},\varphi)$ and 
$D(y,\sigma)=E(p\tilde{\alpha}(y)k,\sigma+\frac{1}{2},\varphi)-
E_P(p\tilde{\alpha}(y)k,\sigma+\frac{1}{2},\varphi).$  Let $a=4$, 
$b=4$, $c=-4$ and $d=0.$  The main analytic lemma applies.  
The asymptotic formula it produces is

$$-\beta \sim \frac{( \varphi,\varphi)}
{2 L(2,\varphi \times \varphi)}
\left|\frac{\det h_1}{\det h_2}\right|^{-2}
y^{-4},$$
where $(\varphi,\varphi)$ is the Petersson inner product of 
$\varphi$ with itself.
\subsection{Example Two: $GL(3)$ over $\Q$}

Now let us consider the case $k=\Q$, and $\G=GL(3,\A)$.  We take $P=P_0$ to
be the Borel subgroup of upper triangular matrices and $M=T_0$ to be the
torus of diagonal matrices.  The associated set of simple positive roots
is $\{\alpha,\theta\}$, defined by
$$\alpha\left(\left(\begin{smallmatrix}t_1&&\\&t_2&\\&&t_3\end{smallmatrix}\right)\right)
=\frac{t_1}{t_2}\hskip .5in
\theta\left(\left(\begin{smallmatrix}t_1&&\\&t_2&\\&&t_3\end{smallmatrix}\right)\right)
=\frac{t_2}{t_3}.
$$
We take $\pi$ to be the trivial representation, and $\phi_\pi$ to be the
constant function 1.  We parametrize 
$\re X_M^{\G}$, by associating the map
$$\I_{\lambda}\left[\left(\begin{smallmatrix}t_1&&\\&t_2&\\&&t_3\end{smallmatrix}\right)\right]=|t_1|^{\lambda_1+1}|t_2|^{\lambda_2}|t_3|^{\lambda_3-1},$$ to the triple $\lambda=(\lambda_1,\lambda_2,\lambda_3) \in \R^3$
such that 
$\lambda_1+\lambda_2+\lambda_3=0.$
 
We work with $\alpha$, and define the map $\tilde{\alpha}$ by
$\tilde{\alpha}(y)=
\left(\begin{smallmatrix}y&&\\&1&\\&&1\end{smallmatrix}\right),$ embedded
at the infinite place.
The Weyl group $W$ of $G$ may be identified with the group of permutation 
matrices.  Then 
$$
W(M,M_{\alpha})=\left\{I,
\left(\begin{smallmatrix}&1&\\1&&\\&&1\end{smallmatrix}\right)
\left(\begin{smallmatrix}&&1\\1&&\\&1&\end{smallmatrix}\right)
\right\}.$$
We compute the intertwining operators as usual, by reducing from the 
general case to the relative rank one case, i.e., to a Levi isomorphic to $GL(2).$  As on $GL(2)$, each root contributes a ratio of Riemann
zeta functions.  Specifically, if $\zeta^*(s)=\pi^{-\frac{s}{2}}\Gamma(\frac{s}{2})\zeta(s)$ is the ``completed'' Riemann zeta function, then we have
\begin{eqnarray*}
M\left(
\left(\begin{smallmatrix}&1&\\1&&\\&&1\end{smallmatrix}\right),
\I_{\lambda}\right)\I_{\lambda}&=&
\frac{\zeta^*(\lambda_1-\lambda_2)}{\zeta^*(\lambda_1-\lambda_2+1)}
\I_{(\lambda_2,\lambda_1,\lambda_3)} \\
M\left(
\left(\begin{smallmatrix}&&1\\1&&\\&1&\end{smallmatrix}\right),
\I_{\lambda}\right)\I_{\lambda}&=&\frac{\zeta^*(\lambda_2-\lambda_3)}{\zeta^*(\lambda_2-\lambda_3+1)}\frac{\zeta^*(\lambda_1-\lambda_3)}{\zeta^*(\lambda_1-\lambda_3+1)}
\I_{(\lambda_3,\lambda_1,\lambda_2)}
\end{eqnarray*}
The value of $E^{M_{\alpha}}$ at a point
$$g=\left(\begin{smallmatrix}1&x_2&x_3\\&1&x_1\\&&1\end{smallmatrix}\right)\left(\begin{smallmatrix}y_1y_2&&\\&y_1&\\&&1\end{smallmatrix}\right) k
\in GL(3,\R)$$
may be given in terms of the Eisenstein series 
$$e(\tau,s)=\sum_{{(c,d)=1}\atop{c>0}}\frac{y^s}{|c\tau+d|^{2s}}$$
for $\tau=x+iy$ 
in the upper half
plane.  Specifically, 
if we let $\tau_1=x_1+iy_1$, then we find that
$$
E^{M_{\alpha}}(\I_{\lambda},\I_{\lambda}\otimes 1)(g)= 
(y_2\sqrt{y_1})^{\lambda_1+1} e(\tau_1,
\frac{\lambda_2-\lambda_3+1}{2}).
$$
It follows that the value of $E_{P_{\alpha}}(\I_{\lambda})$ is given by 
\begin{eqnarray*}
&&(y_2\sqrt{y_1})^{\lambda_1+1} e(\tau_1,
\frac{\lambda_2-\lambda_3+1}{2})\\&&
+(y_2\sqrt{y_1})^{\lambda_2+1} \frac{\zeta^*(\lambda_1-\lambda_2)}{\zeta^*(\lambda_1-\lambda_2+1)}e(\tau_1,\frac{\lambda_1-\lambda_3+1}{2})\\&&
+(y_2\sqrt{y_1})^{\lambda_3+1} \frac{\zeta^*(\lambda_2-\lambda_3)}{\zeta^*(\lambda_2-\lambda_3+1)}
\frac{\zeta^*(\lambda_1-\lambda_3)}{\zeta^*(\lambda_1-\lambda_3+1)}e(\tau_1,\frac{\lambda_1-\lambda_2+1}{2}).\end{eqnarray*}
If we multiply $E(\I_{\lambda},\I_{\lambda}\otimes 1)(g)$ by
$$\zeta^*(\lambda_1-\lambda_2+1)\zeta^*(\lambda_2-\lambda_3+1)\zeta^*(\lambda_1-\lambda_3+1)$$ the resulting
function is essentially the function $G_{\nu_1,\nu_2}$ appearing in 
Bump \cite{yellowbump}.  It is invariant under the six permutations of $\lambda_1,\lambda_2$, and $\lambda_3.$
It's poles are along the lines $\lambda_i-\lambda_j=1,$ and are permuted 
transitively by 
its functional equations.  Let us therefore restrict our attention to the
plane $\lambda_2-\lambda_3=1.$
Now, it's clear that of the three terms that make up $E_{P_{\alpha}}$, the 
first and third will be singular along this hyperplane, while the second will
not.  It is thus clear that
$\{\lambda \mbox{ generic for }\alpha:\wm(\lambda,\alpha)
=\left(\begin{smallmatrix}&1&\\1&&\\&&1
\end{smallmatrix}\right)\}$ 
 is the union of
$\Omega_1:=\{\lambda\in H: \lambda_2>\lambda_1>\lambda_3\}$ and
$\Omega_2:=\{\lambda\in H: \lambda_1<\lambda_3,\lambda_3\neq \lambda_1+1\}$. 
The function $e(\tau,s)$ does not vanish identically in $s$ for any $\tau,$
and its residue is never zero.  Thus, we may use any $p,k.$  
We fix 
$\lambda^0 \in \Omega_1\cup \Omega_2,$ such that 
$e(\tau_1,\frac{\lambda_1^0-\lambda_3^0+1}{2})\neq 0$, and 
 $\lambda^1$ in the hyperplane $\lambda_2-\lambda_3=0$ and define $\Lambda$
by $\Lambda(\sigma)=(1+\sigma)\lambda^0-\sigma \lambda^1.$  There is 
no real loss of generality, since an arbitrary $\Lambda$ may be obtained 
from this one by a simple scaling.  But this ensures that
$\Lambda(\sigma)_2-\Lambda(\sigma)_3=\sigma+1,$ so that the residue of
$\zeta^*(\Lambda(\sigma)_2-\Lambda(\sigma)_3)$ at zero is 1, while that of
$e(\tau, \frac{\Lambda(\sigma)_2-\Lambda(\sigma)_3+1}{2})$ is 
$\frac{1}{\zeta^*(2)}=\frac{6}{\pi}$.
If $\lambda^0\in \Omega_1$,  then
$\wms=I$.  In this case, the asymptotic formula
reads
$$\beta \sim \frac{6\zeta^*(\lambda_1^0-\lambda_2^0+1)}
{\pi \zeta^*(\lambda_1^0-\lambda_2^0)e(\tau_1,\frac{\lambda_1^0-\lambda_3^0+1}{2})(y_2\sqrt{y_1})^{\lambda_2^0-\lambda_1^0}}.$$
On the other hand, if $\lambda^0\in \Omega_\alpha^2$, then 
$\wms=\left(\begin{smallmatrix}&&1\\1&&\\&1&\end{smallmatrix}\right),$
and the asymptotic formula is 
$$\beta \sim \frac{6\zeta^*(\lambda_1^0-\lambda_2^0+1)\zeta^*(\lambda_1^0-\lambda_3^0)e(\tau_1,\frac{\lambda_1^0-\lambda_2^0+1}{2})}
{\pi\zeta^*(\lambda_1^0-\lambda_2^0)\zeta^*(\lambda_1^0-\lambda_3^0+1)e(\tau_1,\frac{\lambda_1^0-\lambda_3^0+1}{2})(y_2\sqrt{y_1})}.$$

\end{document}